\documentclass[smallextended]{svjour3_arxiv}%
\usepackage{graphicx}
\usepackage{amsmath}
\usepackage{amsfonts}
\usepackage{amssymb}
\usepackage{multirow}
\usepackage{tikz}
\usetikzlibrary{arrows.meta}
\usepackage{url}
\usepackage[pdftex]{thumbpdf}
\usepackage[pdftex,pdfstartview=FitH,bookmarks=true,bookmarksnumbered=true,bookmarksopen=true,hypertexnames=false,breaklinks=true,colorlinks=true,linkcolor=blue,anchorcolor=blue,citecolor=blue,filecolor=blue,menucolor=blue]%
{hyperref}%
\providecommand{\U}[1]{\protect\rule{.1in}{.1in}}
\RequirePackage{fix-cm}
\smartqed
\providecommand{\U}[1]{\protect\rule{.1in}{.1in}}
\providecommand{\U}[1]{\protect\rule{.1in}{.1in}}
\newtheorem{algorithm}[theorem]{Algorithm}
\newtheorem{assumption}[theorem]{Assumption}
\begin{document}
\title{Inexact and primal multilevel FETI-DP methods
\thanks
{The work was supported by the U.S. National Science Foundation under award DMS1913201.}
}
\subtitle{a multilevel extension and interplay with BDDC}
\author{Bed\v{r}ich Soused\'{\i}k}
\institute{Bed\v{r}ich Soused\'{\i}k \at
Department of Mathematics and Statistics, University of Maryland, Baltimore County, \\
1000~Hilltop Circle, Baltimore, MD~21250, USA \\
Corresponding author, Tel.: +1-410-455-3298, Fax: +1-410-455-1066\\
\email{sousedik@umbc.edu}}
\date{}
\maketitle\begin{abstract}
We study a framework that allows to solve the coarse problem in the FETI-DP method approximately.
It is based on the saddle-point formulation of the FETI-DP system with a block-triangular preconditioner.
One of the blocks approximates the coarse problem,
for which we use the multilevel BDDC method as the main tool.
This strategy then naturally leads to a version of multilevel FETI-DP method, and
we show that the spectra of the multilevel FETI-DP and BDDC preconditioned operators
are essentially the same. The theory is illustrated by a set of numerical experiments,
and we also present a few experiments when the coarse solve is approximated by algebraic multigrid.
\end{abstract}

\section{Introduction}

\label{sec:introduction}The last two decades in the development of iterative
substructuring class of the domain decomposition methods has been marked by a
significant progress of the two most advanced methods:\ \emph{Finite Element
Tearing and Interconnecting - Dual, Primal} (FETI-DP) and\ \emph{Balancing
Domain Decomposition by Constraints} (BDDC). The FETI-DP method by Farhat et
al.~\cite{Farhat-2001-FDP,Farhat-2000-SDP} was developed from an earlier
version, called FETI (or FETI-1) by Farhat and
Roux~\cite{Farhat-1994-OCP,Farhat-1991-MFE}. The methods from the FETI\ family
are \emph{dual} because they iterate on a system of Lagrange multipliers that
enforce continuity on the interfaces. However, in the FETI-DP some of the
\textquotedblleft coarse\textquotedblright\ variables are treated as
\emph{primal}. The BDDC method was developed by
Dohrmann~\cite{Dohrmann-2003-PSC,Mandel-2003-CBD}, as a fully \emph{primal}
counterpart of the FETI-DP, from an earlier\ BDD\ method by
Mandel~\cite{Mandel-1993-BDD}. Interestingly, two methods that turned out to
be equivalent to BDDC\ were developed independently by
Cros~\cite{Cros-2003-PSC} and by Fragakis and
Papadrakakis~\cite{Fragakis-2003-MHP},
see~\cite{Mandel-2007-BFM,Sousedik-2008-EPD} for a proof. The methods are
closely related, and the essential role here is played by the coarse problem.
The proof that the eigenvalues of the FETI-DP\ and BDDC\ are the same except
possibly for the multiplicity of eigenvalue equal to one, provided that the
coarse components of both methods are the same, was obtained by Mandel et
al.~\cite{Mandel-2005-ATP}, and simpler proofs were derived soon afterwards by
Li and Widlund~\cite{Li-2006-FBB}, and also by Brenner and
Sung~\cite{Brenner-2007-BFW}.

Even though both methods were originally developed for elliptic problems, they
were later extended to other models see, e.g., the
monographs~\cite{Kruis-2006-DDM,Pechstein-2013-FBE,Toselli-2005-DDM} for an
overview. From our current perspective, we find the most important extensions
beyond the original two-level algorithms. The possibility for a multilevel
extension of the BDDC was already mentioned by
Dohrmann~\cite{Dohrmann-2003-PSC}. The theory of three-level BDDC was
developed by Tu~\cite{Tu-2007-TBT3D,Tu-2007-TBT} and extended into a general
multispace BDDC framework (with multilevel BDDC\ as a particular instance)\ by
Mandel et al.~\cite{Mandel-2008-MMB}. For further developments of multilevel
BDDC\ see,
e.g.,~\cite{Kim-2009-TLB,Peng-2018-ABA,Sistek-2013-PIM,Sousedik-2013-NBS,Sousedik-2013-AMB,Tu-2011-TBA,Zampini-2017-MBD}%
. Extensions of the FETI-DP\ are less straightforward because the inverse of
the coarse problem is built directly into the system matrix, and to the best
of our knowledge the literature on this topic is limited. For example, a
hybrid approach to three-level FETI was proposed by Klawonn and
Rheinbach~\cite{Klawonn-2008-HA3}, and a recursive application of the FETI-DP
was recently proposed by Toivanen et al.~\cite{Toivanen-2018-MFM}, see also
Hor\'{a}k~\cite[Chapter~11]{Horak-2007-FBD}.

The FETI-DP\ method is, due to the use of dual variables, suitable in
particular for numerical solution of certain classes of coercive variational
inequalities such as the contact problems see,
e.g.,~\cite{Dostal-2009-OQP,Dostal-2009-SFA,Jarosova-2012-PPT}. However,
factorization of \textcolor{black}{the} coarse problem may become a computational bottleneck when
the number of substructures is large. Our main goal here is to develop a
variant of the FETI-DP method that would allow for an approximate solution of
the coarse problem using a multilevel strategy, but we do not apply the method
recursively as~\cite{Toivanen-2018-MFM}. Instead, we keep the coarse variables
primal. Klawonn and Rheinbach~\cite{Klawonn-2007-IFM} studied inexact FETI-DP
methods based on the saddle-point formulation of the FETI-DP which allows to
move (the approximation of) the inverse of the coarse problem into the
preconditioner. We also use the saddle-point formulation, and the
preconditioner based on an approximation of the block inverse of the system
matrix. However, while~\cite{Klawonn-2007-IFM} applied algebraic multigrid
(AMG) to the coarse problem, our main focus here is the multilevel BDDC and
also our formulation of the algorithms is different since we do not rely on
partial subassembly and a change of variables. This strategy then naturally
leads to a version of multilevel FETI-DP method, which we refer to as
\emph{primal}. We show that the spectra of the preconditioned multilevel
FETI-DP and BDDC\ operators are essentially the same, in the same way as in
the two-level case, and we present numerical experiments that confirm the
theory. We note that an alternative point of view is that our method is a
variant of the multilevel BDDC that allows to use the dual, FETI-DP
formulation\ on the finite element level.

The paper is organized as follows. In Section~\ref{sec:model} we introduce the
model problem and substructuring components, in Section~\ref{sec:methods} we
recall the FETI-DP and BDDC\ methods, in Section~\ref{sec:bddc-dp} we
formulate the inexact and multilevel FETI-DP\ method, in
Section~\ref{sec:analysis} we relate it to the multilevel BDDC, and in
Section~\ref{sec:numerical} we present results of numerical experiments.

\section{Model problem and substructuring components}

\label{sec:model}We introduce a model problem and some standard substructuring
concepts. So, let us assume we are interested in solving a system of
linear\ equations
\begin{equation}
Au=f, \label{eq:Au=f}%
\end{equation}
where $u\in U$, and $A$ is a symmetric, positive definite matrix obtained by a
discretization of an elliptic partial differential equation subject to
appropriate boundary conditions in a bounded domain $\Omega\subset%
\mathbb{R}
^{d}$, $d=2$ or $3$.
Let $\Omega$ be decomposed into $N$ nonoverlapping subdomains$~\Omega^{s}$,
$s=1,\dots,N$, equally called substructures.\ Let each substructure be
obtained as a union of \textcolor{black}{Lagrangian} P1 or Q1 finite elements with mesh size$~h$,
and let$~H$ denote the characteristic size of a subdomain. The first step in
substructuring is elimination of degrees of freedom interior to each
substrure, $u_{I}\in U_{I}$, by so-called \emph{static condensation}, and the
problem is reduced to substructure interfaces. To this end, we denote by
$W_{i}$ the space of degrees of freedom on the substructure$~i$\ interface, by
$w_{i}$ the vector of substructure interface degrees of freedom, and let
\[
W=W^{1}\times\cdots W^{N}.
\]
Next, we denote by $S^{i}:W^{i}\rightarrow W^{i}$ the Schur complement of
subdomain$~i$ after eliminating the interior degrees of freedom, and by
$g^{i}$ the substructure force vector. Also, we introduce $U_{\Gamma}$ as the
space of the global interface degrees of freedom, and denote by $\widehat{W}$
the space of degrees of freedom continuous across the substructure interfaces.
Finally, we introduce~$R^{i}:U_{\Gamma}\rightarrow W^{i}$ as the restriction
map (a zero-one matrix) of global to local degrees of freedom. Then, we
denote
\begin{equation}
S=\left[
\begin{array}
[c]{ccc}%
S^{1} &  & \\
& \ddots & \\
&  & S^{N}%
\end{array}
\right]  \quad w=\left[
\begin{array}
[c]{c}%
w^{1}\\
\vdots\\
w^{N}%
\end{array}
\right]  ,\quad g=\left[
\begin{array}
[c]{c}%
g^{1}\\
\vdots\\
g^{N}%
\end{array}
\right]  ,\quad R=\left[
\begin{array}
[c]{c}%
R^{1}\\
\vdots\\
R^{N}%
\end{array}
\right]  , \label{eq:blocks}%
\end{equation}
where $R$ has a full column rank such that $R^{i}R^{iT}=I$ , and
$\widehat{W}=\operatorname*{range}R$. Then, problem~(\ref{eq:Au=f}) reduced to
interfaces may be written as
\begin{equation}
\widehat{S}w=f_{\Gamma},\qquad\widehat{S}=R^{T}SR,\quad f_{\Gamma}=R^{T}g.
\label{eq:Sw=f}%
\end{equation}
Next ingredient is the averaging operator $E=R^{T}D_{P}$, where $D_{P}%
:W\rightarrow W$ is a given weight matrix such that $ER=I.$ Also, let
$B:W\rightarrow\Lambda$ denote the matrix enforcing the condition
$w\in\widehat{W}$ by $Bw=0$, that is $\operatorname*{range}%
R=\operatorname*{null}B$. In agreement with the block
structure~(\ref{eq:blocks}), we consider
\[
B=\left[
\begin{array}
[c]{ccc}%
B^{1} & \cdots & B^{N}%
\end{array}
\right]  ,\qquad B^{i}=W^{i}\rightarrow\Lambda,\quad\operatorname*{range}%
B=\Lambda.
\]
Let $B_{D}^{T}$ denote a generalized inverse of $B$, that is $B_{D}%
^{T}:\Lambda\rightarrow W$ and $BB_{D}^{T}=I_{\Lambda}$. In computations, the
matrix $B_{D}$ is constructed from $B$\ and $D_{P}$,
see~\cite{Klawonn-2002-DPF,Mandel-2005-ATP}, and, in particular, we will
assume that
\[
B_{D}^{T}B+RE=I.
\]

The FETI-DP and BDDC\ methods are characterized by a selection of \emph{coarse
degrees of freedom }(resp. \emph{constraints}). These can be values at corners
or averages over edges or faces. Let $\widetilde{W}\subset W$ be the subspace
of vectors with continuous coarse degrees of freedom across substructure
interface. Specifically, suppose we are given a space$~X$ and a matrix
$C:W\rightarrow X$, and define
\begin{equation}
\widetilde{W}=\left\{  w\in W:C\left(  I-RE\right)  w=0\right\}  ,
\label{eq:W-tilde-1}%
\end{equation}
where each row of $C$ defines one constraint. The values $Cw$ are called local
coarse degrees of freedom and by~(\ref{eq:W-tilde-1}) they have zero jumps on
adjacent substructures. To represent the global coarse degrees of freedom, we
use a space$~U_{c}$ and a matrix$~R_{c}$ with full column rank, such that
\[
\widetilde{W}=\left\{  w\in W:\exists u_{c}\in U_{c}:Cw=R_{c}u_{c}\right\}  .
\]
Next, let $\Phi$ be a matrix that defines the basis for the coarse problem,
given as
\[
\Phi=\left[
\begin{array}
[c]{c}%
\Phi^{1}\\
\vdots\\
\Phi^{N}%
\end{array}
\right]  ,\qquad C\Phi=R_{c}.
\]
\textcolor{black}{Here we use energy minimal coarse basis functions$~\Psi$ that minimize 
energy of the form $x^T \Phi^{T}S\Phi x$ subject to constraints $Cx=e$ for some vector $e$ given for each basis function as a column of matrix $R_c$.
Specifically, the basis functions are found by solving saddle-point problems}
\begin{equation}
\left[
\begin{array}
[c]{cc}%
S & C^{T}\\
C & 0
\end{array}
\right]  \left[
\begin{array}
[c]{c}%
\Psi\\
\Lambda
\end{array}
\right]  =\left[
\begin{array}
[c]{c}%
0\\
R_{c}%
\end{array}
\right]  . \label{eq:energy-minimization}%
\end{equation}
In implementation, we use a decomposition of the space $\widetilde{W}$,
see~\cite[Lemma~8]{Mandel-2005-ATP},
\begin{equation}
\widetilde{W}=\widetilde{W}_{\Pi}\oplus\widetilde{W}_{\Delta},
\label{eq:dharm-dec-interface}%
\end{equation}
where $\widetilde{W}_{\Pi}=\operatorname*{range}\Phi$ is the \emph{primal}
space, which gives rise to the (global) coarse problem, and $\widetilde{W}%
_{\Delta}=\operatorname*{null}C$ is the \emph{dual} space with coarse degrees
of freedom equal to zero, which gives rise to independent subdomain problems.

\section{FETI-DP and BDDC\ methods}

\label{sec:methods}

\subsection{FETI-DP}

Problem~(\ref{eq:Sw=f}) may be equivalently\ written as
\begin{equation}
\min_{w\in\widetilde{W}}\max_{\lambda\in\Lambda}\mathcal{L}(w,\lambda),
\label{eq:nested-minimization}%
\end{equation}
where $\mathcal{L}\left(  w,\lambda\right)  $ is the \textcolor{black}{Lagrangian} 
\[
\mathcal{L}\left(  w,\lambda\right)  =\frac{1}{2}w^{T}Sw-w^{T}E^{T}f_{\Gamma
}+w^{T}B^{T}\lambda.
\]
Problem~(\ref{eq:nested-minimization}) is further equivalent to the dual
problem $\max_{\lambda\in\Lambda}\mathcal{F}(\lambda)$, where
\begin{equation}
\mathcal{F}(\lambda)=\min_{w\in\widetilde{W}}\mathcal{L}(w,\lambda)
\label{eq:dual-functional}%
\end{equation}
is the dual functional.

Next, we use the splitting~(\ref{eq:dharm-dec-interface}) and introduce
Lagrange multipliers$~\mu$ to enforce the zero values of coarse degrees of
freedom in$~\widetilde{W}_{\Delta}$. Then,
functional~(\ref{eq:dual-functional}) can be written with minimization
over$~\widetilde{W}$ instead of$~\widetilde{W}_{\Delta}$ as
\begin{equation}
\mathcal{F}(\lambda)=\min_{w_{\Delta}\in\widetilde{W}}\min_{w_{\Pi}%
\in\widetilde{W}_{\Pi}}\sup_{\mu}\mathcal{L}(w_{\Delta}+w_{\Pi},\lambda
)+w_{\Delta}^{T}C^{T}\mu.\label{eq:augmented-nested-minimization}%
\end{equation}
Writing $w_{\Pi}=\Phi u_{c}$,
and differentiating with respect to $w_{\Delta}$, $\mu$, $u_{c}$, $\lambda$,
we get a linear system, with a solution equvalent to solving the dual problem,
which is
\begin{equation}
\left[
\begin{array}
[c]{cccc}%
S & C^{T} & S\Phi & B^{T}\\
C & 0 & 0 & 0\\
\Phi^{T}S^{T} & 0 & \Phi^{T}S\Phi & \Phi^{T}B^{T}\\
B & 0 & B\Phi & 0
\end{array}
\right]  \left[
\begin{array}
[c]{c}%
w_{\Delta}\\
\mu\\
u_{c}\\
\lambda
\end{array}
\right]  =\left[
\begin{array}
[c]{c}%
E^{T}f_{\Gamma}\\
0\\
\Phi^{T}E^{T}f_{\Gamma}\\
0
\end{array}
\right]  .\label{eq:feti-dp-full}%
\end{equation}
Eliminating $w_{\Delta}$ and $\mu$ from the first two equations, we derive the
so-called \emph{saddle-point }FETI-DP formulation, which can be written as
\begin{equation}
\left[
\begin{array}
[c]{cc}%
S_{c} & \Psi^{T}B^{T}\\
B\Psi & -BS_{\Delta}^{-1}B^{T}%
\end{array}
\right]  \left[
\begin{array}
[c]{c}%
u_{c}\\
\lambda
\end{array}
\right]  =\left[
\begin{array}
[c]{c}%
\Psi^{T}E^{T}f_{\Gamma}\\
-BS_{\Delta}^{-1}E^{T}f_{\Gamma}%
\end{array}
\right]  ,\label{eq:bddc-dp}%
\end{equation}
where $S_{\Delta}^{-1}$ corresponds to independent subdomain solves
\begin{equation}
S_{\Delta}^{-1}=\left[
\begin{array}
[c]{c}%
I\\
0
\end{array}
\right]  ^{T}\left[
\begin{array}
[c]{cc}%
S & C^{T}\\
C & 0
\end{array}
\right]  ^{-1}\left[
\begin{array}
[c]{c}%
I\\
0
\end{array}
\right]  ,\label{eq:sub_solve}%
\end{equation}
so that $\Psi=\left(  I-S_{\Delta}^{-1}S\right)  \Phi$ are the energy
minimizing coarse basis functions, the same as
in~(\ref{eq:energy-minimization}), cf.~\cite[eqs.~($22$) and ($25$%
)]{Mandel-2005-ATP}, and $S_{c}$ is the coarse matrix
\begin{equation}
S_{c}=\Phi^{T}S\Phi-\left[
\begin{array}
[c]{c}%
S\Phi\\
0
\end{array}
\right]  ^{T}\left[
\begin{array}
[c]{cc}%
S & C^{T}\\
C & 0
\end{array}
\right]  ^{-1}\left[
\begin{array}
[c]{c}%
S\Phi\\
0
\end{array}
\right]  =\Psi^{T}S\Psi,\label{eq:coarse_matrix}%
\end{equation}
where the second equality follows from~\cite[Theorem~$3$]{Mandel-2005-ATP}.

We note that eliminating$~u_{c}$ from (\ref{eq:bddc-dp}) yields the
\emph{original} FETI-DP\ system
\begin{equation}
B\left(  \Psi S_{c}^{-1}\Psi^{T}+S_{\Delta}^{-1}\right)  B^{T}\lambda=B\left(
\Psi S_{c}^{-1}\Psi^{T}+S_{\Delta}^{-1}\right)  E^{T}f_{\Gamma}%
.\label{eq:feti-dp-orig}%
\end{equation}
The original FETI-DP is the method of preconditioned conjugate gradients
applied to~(\ref{eq:feti-dp-orig}), used typically with the so-called
Dirichlet preconditioner
\begin{equation}
M_{D}=B_{D}SB_{D}^{T}.\label{eq:Dirichlet-preconditioner}%
\end{equation}

\subsection{BDDC}

The BDDC preconditioner was originally formulated by
Dohrmann~\cite{Dohrmann-2003-PSC} for problem~(\ref{eq:Au=f}). Formulations
for the reduced problem~(\ref{eq:Sw=f}) followed
in~\cite{Mandel-2003-CBD,Mandel-2005-ATP}. The BDDC is the method of
preconditioned conjugate gradients, in which the (two-level)
BDDC\ preconditioner for problem~(\ref{eq:Sw=f}) is given as,
cf.~\cite[Lemma~27]{Mandel-2005-ATP},
\begin{equation}
M_{\widehat{S}}=EHE^{T}, \label{eq:bddc}%
\end{equation}
where%
\begin{equation}
H=\left[
\begin{array}
[c]{c}%
\Psi^{T}\\
I
\end{array}
\right]  ^{T}\left[
\begin{array}
[c]{cc}%
S_{c} & 0\\
0 & S_{\Delta}%
\end{array}
\right]  ^{-1}\left[
\begin{array}
[c]{c}%
\Psi^{T}\\
I
\end{array}
\right]  =\Psi S_{c}^{-1}\Psi^{T}+S_{\Delta}^{-1}. \label{eq:H}%
\end{equation}
The preconditioner formulated for problem~(\ref{eq:Au=f}) includes in addition
a pre-processing step with a static condensation of the residual and
post-processing step with recovery of the approximate solution corresponding
to the subdomain interiors. Both formulations are equivalent, and in
particular by~\cite[Theorem~14]{Mandel-2008-MMB}, the eigenvalues of the two
preconditioned operators are the same except possibly for multiplicity of
eigenvalue equal to one.

\subsubsection{Multilevel BDDC}

In case when the number of subdomains in large, solving the \emph{coarse
problem}\ becomes a bottleneck. This motivated the introduction of the
multilevel extensions~\cite{Mandel-2008-MMB,Tu-2007-TBT3D,Tu-2007-TBT} that
solve the coarse problem only approximately by applying the two-level BDDC
recursively. Since we will apply this idea in the FETI-DP\ framework, we
briefly recall the Multilevel BDDC from~\cite{Mandel-2008-MMB}.
The substructuring components\ will be denoted by an additional
subscript$~_{1}$ as $U_{1}$, $\widetilde{W}_{\Pi1}$\ and $\Omega_{1}^{s}$,
$s=1,\dots,N_{1}$, etc., and called level$~1$. Level$~1$ coarse problem will
be called level$~2$ problem. Because it has an analogous structure as the
original problem on level$~1$, we put $U_{2}=\widetilde{W}_{\Pi1}$, and we
introduce the substructuring components for level$~2$ in the same way as we
introduced level$~1$ components.
Generally, in a design of $L$-level method, we repeat this process recursively
for levels $\ell=1,\dots,L-1$. On each decomposition level$~\ell$, we assume
that the subdomains$~\Omega_{\ell}^{s}$, $s=1,\dots,N_{\ell}$, have
characteristic size $H_{\ell}$ and form a conforming triangulation of the
domain$~\Omega$. \ Level$~\ell-1$ substructures become level$~\ell$ elements,
level$~\ell-1$ coarse degrees of freedom become level$~\ell$ degrees of
freedom, and $\Gamma_{\ell}\subset\Gamma_{\ell-1}$. The finite element level
will be denoted as level$~0$\ and $H_{0}=h$. 
\textcolor{black}{An example of a decomposition is shown in Figure~\ref{fig:decomposition}.}

\begin{figure}[tbh]
\begin{center}%
\begin{tikzpicture}[scale=.14] 
\foreach \i in {0,...,9} { \draw [very thin,gray] (\i,0) -- (\i,9);} \foreach \i in {0,...,9} { \draw [very thin,gray] (0,\i) -- (9,\i); } 
\foreach \i in {0,...,9} { \draw [thick,red] (3*\i,0) -- (3*\i,18);} \foreach \i in {0,...,6} { \draw [thick,red] (0,3*\i) -- (27,3*\i); } 
\foreach \i in {0,...,3} { \draw [very thick,blue] (9*\i,0) -- (9*\i,27);} \foreach \i in {0,...,3} { \draw [very thick,blue] (0,9*\i) -- (27,9*\i); } 
\draw (-21,0) -- (0,0); \draw (-2.8,1) -- (0,1);  \draw [->] (-2,-1.5) -- (-2,0); \draw (-2,0) -- (-2,1); \draw [->] (-2,2.5) -- (-2,1); \node at (-6.5,1.5) {$H_0=h$};
\draw (-15,3) -- (0,3); \draw [->] (-14,-2) -- (-14,0); \draw (-14,0) -- (-14,3); \draw [->] (-14,5) -- (-14,3); \node at (-16,1.5) {$H_1$};
\draw (-21,9) -- (0,9); \draw [->] (-20,-2) -- (-20,0); \draw (-20,0) -- (-20,9); \draw [->] (-20,12) -- (-20,9); \node at (-22,4.5) {$H_2$};
\draw [{Circle[color=red]}-] (25,1.5) -- (30,1.5);  \node[font = {\large}, red] at (33,1.5) {$\Omega_1^9$}; 
\draw [{Circle[color=blue]}-] (22,22.5) -- (30,22.5);  \node[font = {\large}, blue] at (33,22.5) {$\large \Omega_2^9$}; 
\end{tikzpicture} 
\end{center}
\caption{An example of a uniform decomposition for a three-level method with $H_\ell/H_{\ell-1}=3$, where $\ell=1,2$. 
Level~$1$ subdomains are shown in red and level~$2$ subdomains in blue.}%
\label{fig:decomposition}%
\end{figure}

The shape functions on
level$~\ell$ are determined by minimization of energy with respect to
level$~\ell-1$ shape functions, subject to the value of exactly one
level$~\ell$ degree of freedom being one and others level $\ell$ degrees of
freedom being zero. The minimization is done on each level$~\ell$ element
(level$~\ell-1$ substructure) separately, so the values of level$~\ell-1$
degrees of freedom are in general discontinuous between level$~\ell-1$
substructures, and only the values of level$~\ell$ degrees of freedom between
neighboring level$~\ell$ elements coincide.
The hierarchy of the spaces for $L$-level BDDC\ method can be written as
\[%
\begin{array}
[c]{ccccccccccccccccccc}
&  & U_{I1} &  &  &  &  &  &  &  &  &  &  &  & U_{IL-1} &  &  &  & \\
U_{1} & = & \oplus &  &  &  & \widetilde{W}_{\Pi1} & = & U_{2} &  & \cdots &
& U_{L-1} & = & \oplus &  &  &  & \widetilde{W}_{\Pi L-1}\\
&  & \widehat{W}_{1} &
\genfrac{}{}{0pt}{}{\genfrac{}{}{0pt}{}{E_{1}}{\leftarrow}}{\subset}%
& \widetilde{W_{1}} & = & \oplus &  &  &  &  &  &  &  & \widehat{W}_{L-1} &
\genfrac{}{}{0pt}{}{\genfrac{}{}{0pt}{}{E_{L-1}}{\leftarrow}}{\subset}%
& \widetilde{W}_{L-1} & = & \oplus\\
&  &  &  &  &  & \widetilde{W}_{\Delta1} &  &  &  &  &  &  &  &  &  &  &  &
\widetilde{W}_{\Delta L-1}%
\end{array}
.
\]
The two-level BDDC\ preconditioner for problem~(\ref{eq:Sw=f}) entails finding
corrections in the spaces$~\widetilde{W}_{\Delta1}$ and$~\widetilde{W}_{\Pi1}%
$, and the multilevel BDDC\ for the same problem entails computing corrections
in the spaces as above, that is
\begin{equation}
V_{1}=\widetilde{W}_{\Delta1},\quad V_{2}=U_{I2},\quad V_{3}=\widetilde{W}%
_{\Delta2},\quad\dots\quad V_{2(L-1)}=\widetilde{W}_{\Pi L-1},
\label{eq:ml-spaces}%
\end{equation}
such that
\begin{equation}
\widetilde{W}_{1}=\left(  \widetilde{W}_{\Delta1}\oplus\widetilde{W}_{\Pi
1}\right)  \subset\widetilde{V}=\sum\nolimits_{k=1}^{2(L-1)}V_{k}\subset W.
\label{eq:ml-spaces-embedding}%
\end{equation}
The precise formulation of the multilevel BDDC preconditioner can be found
in~\cite[Algorithm 17]{Mandel-2008-MMB}. More generally, as an analogy
to~(\ref{eq:bddc})--(\ref{eq:H}), we may consider an \emph{inexact}
BDDC\ preconditioner given as
\begin{equation}
\widetilde{M}_{\widehat{S}}=E\widetilde{H}E^{T}, \label{eq:mBDDC}%
\end{equation}
where
\begin{equation}
\widetilde{H}=\left[
\begin{array}
[c]{c}%
\Psi^{T}\\
I
\end{array}
\right]  ^{T}\left[
\begin{array}
[c]{cc}%
M_{c} & 0\\
0 & S_{\Delta}^{-1}%
\end{array}
\right]  \left[
\begin{array}
[c]{c}%
\Psi^{T}\\
I
\end{array}
\right]  =\Psi M_{c}\Psi^{T}+S_{\Delta}^{-1}, \label{eq:H-tilde}%
\end{equation}
and $M_{c}$ is a preconditioner for the coarse matrix$~S_{c}$. The multilevel
BDDC\ applied to the coarse problem is then a particular example of$~M_{c}$.

Finally, we note that the derivation of the FETI-DP~(\ref{eq:feti-dp-orig})
provides an alternative proof to~\cite[Lemma~$27$]{Mandel-2005-ATP}, which
states that the preconditioned BDDC\ and FETI-DP\ operators can be written,
respectively, as
\begin{align}
M_{\widehat{S}}\widehat{S}  &  =\left(  EHE^{T}\right)  \left(  R^{T}%
SR\right)  ,\label{eq:PA}\\
M_{D}F  &  =\left(  B_{D}SB_{D}^{T}\right)  \left(  BHB^{T}\right)  ,
\label{eq:MF}%
\end{align}
which is used in the proof of equivalence of the spectra of the two operators.

\begin{theorem}
[{\cite[Theorem~$26$]{Mandel-2005-ATP}}]
\label{thm:equiv}
The eigenvalues of the preconditioned
operators of the BDDC and FETI-DP\ methods given, respectively, by~(\ref{eq:PA}) and~(\ref{eq:MF})
are the same except possibly for the eigenvalues equal to zero and one.
\end{theorem}

\begin{proof}
See Appendix.
\end{proof}

Our next goal is to design a multilevel version of the FETI-DP\ method, based
on formulation~(\ref{eq:bddc-dp}), and relate its spectrum to that of
Multilevel BDDC.

\section{Inexact and multilevel FETI-DP methods}

\label{sec:bddc-dp}Since the inverse of the coarse matrix$~S_{c}$ is built
into the original FETI-DP operator in (\ref{eq:feti-dp-orig}), its application
may become a bottleneck when the number of substructures is large. Therefore,
we work with the saddle-point formulation~(\ref{eq:bddc-dp}) that allows to
solve coarse problem inexactly, which is a strategy proposed by Klawonn and
Rheinbach~\cite{Klawonn-2007-IFM}.

The proposed preconditioner is based on the block inverse
\begin{align*}
\left[
\begin{array}
[c]{cc}%
X & Y^{T}\\
Y & -Z
\end{array}
\right]  ^{-1}  &  =\left[
\begin{array}
[c]{cc}%
I & -X^{-1}Y^{T}\\
0 & I
\end{array}
\right]  \left[
\begin{array}
[c]{cc}%
X^{-1} & 0\\
0 & -S_{X}^{-1}%
\end{array}
\right]  \left[
\begin{array}
[c]{cc}%
I & 0\\
-YX^{-1} & I
\end{array}
\right]  =\\
&  =\left[
\begin{array}
[c]{cc}%
I & -X^{-1}Y^{T}\\
0 & I
\end{array}
\right]  \left[
\begin{array}
[c]{cc}%
X^{-1} & 0\\
S_{X}^{-1}YX^{-1} & -S_{X}^{-1}%
\end{array}
\right]  ,
\end{align*}
and $S_{X}=Z+YX^{-1}Y^{T}$. For the inverse of matrix in~(\ref{eq:bddc-dp}),
we then get
\begin{equation}
\left[
\begin{array}
[c]{cc}%
S_{c} & \Psi^{T}B^{T}\\
B\Psi & -BS_{\Delta}^{-1}B^{T}%
\end{array}
\right]  ^{-1}=\left[
\begin{array}
[c]{cc}%
I & -S_{c}^{-1}\Psi^{T}B^{T}\\
0 & I
\end{array}
\right]  \left[
\begin{array}
[c]{cc}%
S_{c}^{-1} & 0\\
F^{-1}B\Psi S_{c}^{-1} & -F^{-1}%
\end{array}
\right]  , \label{eq:block-dec-bddc-dp}%
\end{equation}
and $F=B\left(  \Psi S_{c}^{-1}\Psi^{T}+S_{\Delta}^{-1}\right)  B^{T}$ is the
same as the FETI-DP matrix in~(\ref{eq:feti-dp-orig}). The matrix on the right
in~(\ref{eq:block-dec-bddc-dp}) then motivates the choice of a preconditioner
\begin{equation}
M_{F}=\left[
\begin{array}
[c]{cc}%
M_{c} & 0\\
M_{D}B\Psi M_{c} & -M_{D}%
\end{array}
\right]  , \label{eq:M_D}%
\end{equation}
where$~M_{D}$ is the Dirichlet
preconditioner~(\ref{eq:Dirichlet-preconditioner}), and $M_{c}$ is the
preconditioner for the coarse matrix$~S_{c}$. In general, we may consider any
suitable$~M_{c}$ giving rise to \emph{inexact} FETI-DP. For example, Klawonn
and Rheinbach~\cite{Klawonn-2007-IFM} used algebraic multigrid (AMG). Here,
our main focus is on application of multilevel BDDC in place of$~M_{c}$, which
we then refer to as \emph{primal multilevel} FETI-DP method. In
implementation, the application of the preconditioner~$M_{F}$ entails one
application of both~$M_{c}$ and~$M_{D}$, which may be computed as follows.

\begin{algorithm}
\label{alg:M_D}The preconditioner$~M_{F}:$ $\left(  r_{c},r_{\lambda}\right)
\mapsto\left(  u_{c},\lambda\right)  $
is applied as:\newline First, apply the preconditioner $M_{c}$ as $u_{c}%
=M_{c}r_{c}$, and compute%
\begin{align*}
x  &  =B\Psi u_{c},\\
y  &  =x-r_{\lambda}.
\end{align*}
Then, apply the preconditioner $M_{D}$ as $\lambda=M_{D}y$, and concatenate
the results
\[
\left[
\begin{array}
[c]{c}%
u_{c}\\
\lambda
\end{array}
\right]  .
\]

\end{algorithm}

\begin{remark}
Algorithm~\ref{alg:M_D} is closely related to the second algorithm by
Klawonn and Rheinbach~\cite{Klawonn-2007-IFM}.
However, our formulation of FETI-DP does not use partial
subassembly and a change of variables,
and we consider multilevel BDDC besides algebraic multigrid as a preconditioner~$M_{c}$
for the coarse matrix~$S_c$.
\end{remark}

\section{Analysis of the methods}

\label{sec:analysis}The saddle-point FETI-DP matrix from~(\ref{eq:bddc-dp})
preconditioned by$~M_{F}$\ from~(\ref{eq:M_D})~is%
\begin{equation}
\left[
\begin{array}
[c]{cc}%
M_{c}S_{c} & M_{c}\Psi^{T}B^{T}\\
M_{D}B\Psi\left(  M_{c}S_{c}-I\right)  & M_{D}B\left(  \Psi M_{c}\Psi
^{T}+S_{\Delta}^{-1}\right)  B^{T}%
\end{array}
\right]  . \label{eq:M_D_bddc-dp}%
\end{equation}
First, note that with the exact coarse solve, i.e., setting $M_{c}=S_{c}^{-1}%
$, the matrix in position (1,1) is the identity so the matrix in position
(2,1) vanishes, and the matrix in position (2,2) is the same as the
preconditioned FETI-DP operator~(\ref{eq:MF}). Finally, multiplying the
right-hand side of~(\ref{eq:bddc-dp}) by$~M_{F}$ we may check that solving the
second set of equations therein is the same as solving the original
FETI-DP\ problem~(\ref{eq:feti-dp-orig}) with the Dirichlet
preconditioner~$M_{D}$.

Next, let us observe that the inexact BDDC preconditioner~(\ref{eq:mBDDC})
with $\widetilde{H}$ from~(\ref{eq:H-tilde})\ and the entry in position (2,2)
of matrix (\ref{eq:M_D_bddc-dp}) with the Dirichlet preconditioner$~M_{D}$
defined by~(\ref{eq:Dirichlet-preconditioner}) may be written, respectively,
as
\begin{align}
\widetilde{M}_{\widehat{S}}\widehat{S}  &  =\left(  E\widetilde{H}%
E^{T}\right)  \left(  R^{T}SR\right)  ,\label{eq:mBDDC-prec-op}\\
M_{D}\widetilde{F}  &  =\left(  B_{D}SB_{D}^{T}\right)  \left(  B\widetilde{H}%
B^{T}\right)  . \label{eq:(2,2)}%
\end{align}
That is (\ref{eq:mBDDC-prec-op})--(\ref{eq:(2,2)}) is an analogy to
(\ref{eq:PA})--(\ref{eq:MF}), just with $\widetilde{H}$ in place of $H$.

Let us now focus on the primal multilevel FETI-DP\ method, which combines the
dual (FETI-DP type of) approach to the solve in~$\widetilde{W}_{\Delta}$ and
applies the (primal) multilevel BDDC preconditioner to the coarse solve
in~$\widetilde{W}_{\Pi}$.

\begin{assumption}
\label{assum:coarse-split}We will assume that the spaces $V_{k}$,
$k=1,\dots,2(L-1)$, are constructed so that the associated linear systems are
nonsingular, and the multilevel BDDC uses exact solves in these subspaces, so
that
\begin{equation}
\left.  \widetilde{H}S\right\vert _{\widetilde{V}\subset W}=I.\label{eq:HS}%
\end{equation}

\end{assumption}

For the two-level BDDC method the Assumption~\ref{assum:coarse-split} holds
due to the exact coarse solve, see also~\cite[Lemma~28]{Mandel-2005-ATP}. For
multispace and multilevel BDDC see \cite[Corollary~4, Remark~5 and
Lemma~18]{Mandel-2008-MMB}. 

Then, a variant of Theorem~\ref{thm:equiv} also holds for the two
operators~(\ref{eq:mBDDC-prec-op})--(\ref{eq:(2,2)}).

\begin{lemma}
\label{lem:equiv-mlevel}
Let $\widetilde{H}$ be defined as in~(\ref{eq:H-tilde}%
) and so that assumption~(\ref{eq:HS}) is satisfied. Then the eigenvalues
of the preconditioned operators given by (\ref{eq:mBDDC-prec-op}) and
(\ref{eq:(2,2)}) are the same except possibly for the eigenvalues equal to
zero and one.
\end{lemma}

\begin{proof}
See Appendix.
\end{proof}

\textcolor{black}{Let us now use Lemma~\ref{lem:equiv-mlevel} to examine the eigenvalues of the preconditioned operator~(\ref{eq:M_D_bddc-dp}). 
The corresponding eigenvalue problem can be written~as} 

\begin{equation}%
\begin{array}
[c]{ccccc}%
M_{c}S_{c}u & + & M_{c}\Psi^{T}B^{T}p & = & \lambda u,\\
M_{D}B\Psi\left(  M_{c}S_{c}-I\right)  u & + & M_{D}B\left(  \Psi M_{c}%
\Psi^{T}+S_{\Delta}^{-1}\right)  B^{T}p & = & \lambda p.
\end{array}
\label{eq:1}%
\end{equation}
We begin by considering possible eigenvalues of $M_{c}S_{c}$. There are two
possible cases for $M_{c}S_{c}u=\mu u$ where $(\mu,u)$ is the corresponding
eigenpair: $\mu=1$ or $\mu\neq1$. If $\mu=1$, then $M_{c}S_{c}u=u$ and,
assuming $\lambda\neq1$, we get from the first equation $u=1/(\lambda
-1)M_{c}\Psi^{T}B^{T}p$. Then, the second equation of~(\ref{eq:1}) simplifies
to the eigenvalue problem for~(\ref{eq:(2,2)}), that is $M_{D}\widetilde{F}%
p=\lambda p$, which has essentially the same eigenvalues as the multilevel
BDDC preconditioned operator $\widetilde{M}_{\widehat{S}}\widehat{S}$ by
Lemma~\ref{lem:equiv-mlevel}. If $\mu\neq1$, then $M_{c}S_{c}u=\mu u$ and,
assuming $\lambda\neq\mu$, we get from the first equation $u=1/(\lambda
-\mu)M_{c}\Psi^{T}B^{T}p$. Substituting this into the second equation
of~(\ref{eq:1}), we obtain
\[
M_{D}B\left(  \frac{\mu-1}{\lambda-\mu}\Psi M_{c}\Psi^{T}+\Psi M_{c}\Psi
^{T}+S_{\Delta}^{-1}\right)  B^{T}p=\lambda p,
\]
which reduces with $\lambda=1$ to $M_{D}BS_{\Delta}^{-1}B^{T}p=\lambda p$.
Therefore, Theorem~\ref{thm:equiv} for two-level methods translates to the
multilevel BDDC\ and FETI-DP\ methods.

\begin{theorem}
The eigenvalues of the preconditioned operators of the multilevel BDDC and primal multilevel FETI-DP methods
are the same except possibly for the eigenvalues equal to zero and one.
\end{theorem}

For further discussion of multilevel BDDC\ method we refer
to~\cite{Mandel-2008-MMB}. It is well known that the eigenvalues of the
operators provide only asymptotic insight into the convergence of GMRES, which
is the main iterative method used in the numerical experiments. The bound on
the error at every step involves the condition number of the operator
eigenvectors see, e.g.,~\cite{Elman-2014-FEF,Saad-2003-IMS,Simoncini-2004-BTP}%
, which is difficult to estimate. Nevertheless, the numerical experiments
presented in the next section illustrate that the convergence of GMRES applied
to either multilevel FETI-DP\ or BDDC\ method is quite comparable to that of
the preconditioned conjugate gradients applied to the multilevel BDDC. We also
note that if the preconditioned operator is symmetric, positive definite in
some inner product then it may be also used in conjugate gradient
method~\cite{Klawonn-1998-BTP,Klawonn-2007-IFM}.

\section{Numerical experiments}

\label{sec:numerical}The methods were implemented in \textsc{Matlab}. We used
the right-preconditioned GMRES\ method with no restarts and (for multilevel
BDDC also) preconditioned conjugate gradients, with relative residual
tolerance$~10^{-8}$ as the stopping criterion and zero initial guess. To
estimate the largest eigenvalues of the preconditioned operators, we used the
\textsc{Matlab} function \texttt{eigs}.

First, we presents results of numerical experiments for the model problem
corresponding to the scalar second-order elliptic problem with zero Dirichlet
boundary conditions on a square domain discretized by linear quadrilateral
finite elements. The domain was uniformly divided into substructures with
fixed $H_{\ell}/H_{\ell-1}$ ratio on each level$~\ell$. Table~\ref{tab:c}
shows the largest eigenvalues of the preconditioned operators and the
iteration counts of the multilevel FETI-DP\ and BDDC\ methods set up using
corner constraints and varying numbers of levels$~L$ and coarsening ratios
$H_{\ell}/H_{\ell-1}$, $\ell=1,\dots,L-1$. Table~\ref{tab:c+f} then shows
results for the same set of problems and the two methods set up using corner
constraints combined with arithmetic averages over edges. We used
right-preconditioned GMRES for the FETI-DP\ method~(\ref{eq:bddc-dp})
preconditioned by$~M_{F}$ from~(\ref{eq:M_D}), in which multilevel
BDDC\ preconditioner was applied as$~M_{c}$, and for the multilevel
BDDC\ method we used also preconditioned conjugate gradients (PCG). In
addition, for the FETI-DP\ formulation, we studied performance of a
block-diagonal variant of the preconditioner$~M_{F}$, that is
\begin{equation}
\left[
\begin{array}
[c]{cc}%
M_{c} & 0\\
0 & -M_{D}%
\end{array}
\right]  ,\label{eq:M_BD}%
\end{equation}
which from the theory in~\cite{Murphy-2000-NPI} may be nearly as efficient.
From the tables it can be seen the largest eigenvalues of both multilevel
FETI-DP\ and BDDC\ preconditioned operators were the same. Also, the iteration
counts of GMRES\ were essentially the same for the two methods, and the same
trend can be observed for the PCG\ method. However, in case of the FETI-DP
method\ with the block-diagonal preconditioner the iteration counts of
GMRES\ are higher in all cases. Since the overhead associated with the
application of$~M_{F}$ is small, cf. Algorithm~\ref{alg:M_D}, the
block-diagonal preconditioner appears to be less suitable. 

\begin{figure}[ptbh]
\begin{center}%
\begin{tabular}
[c]{cc}%
\includegraphics[width=5.5cm]{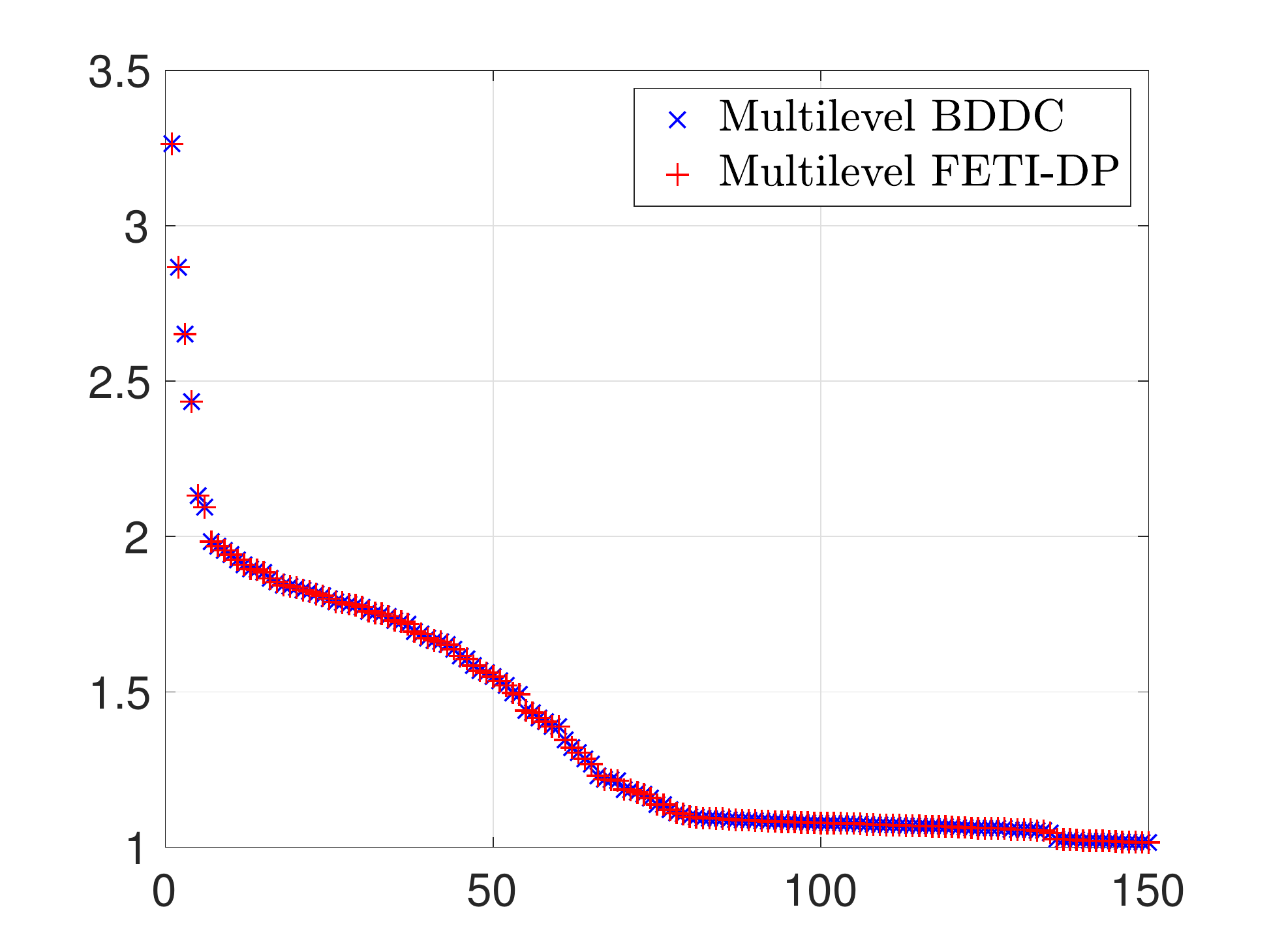} &
\includegraphics[width=5.5cm]{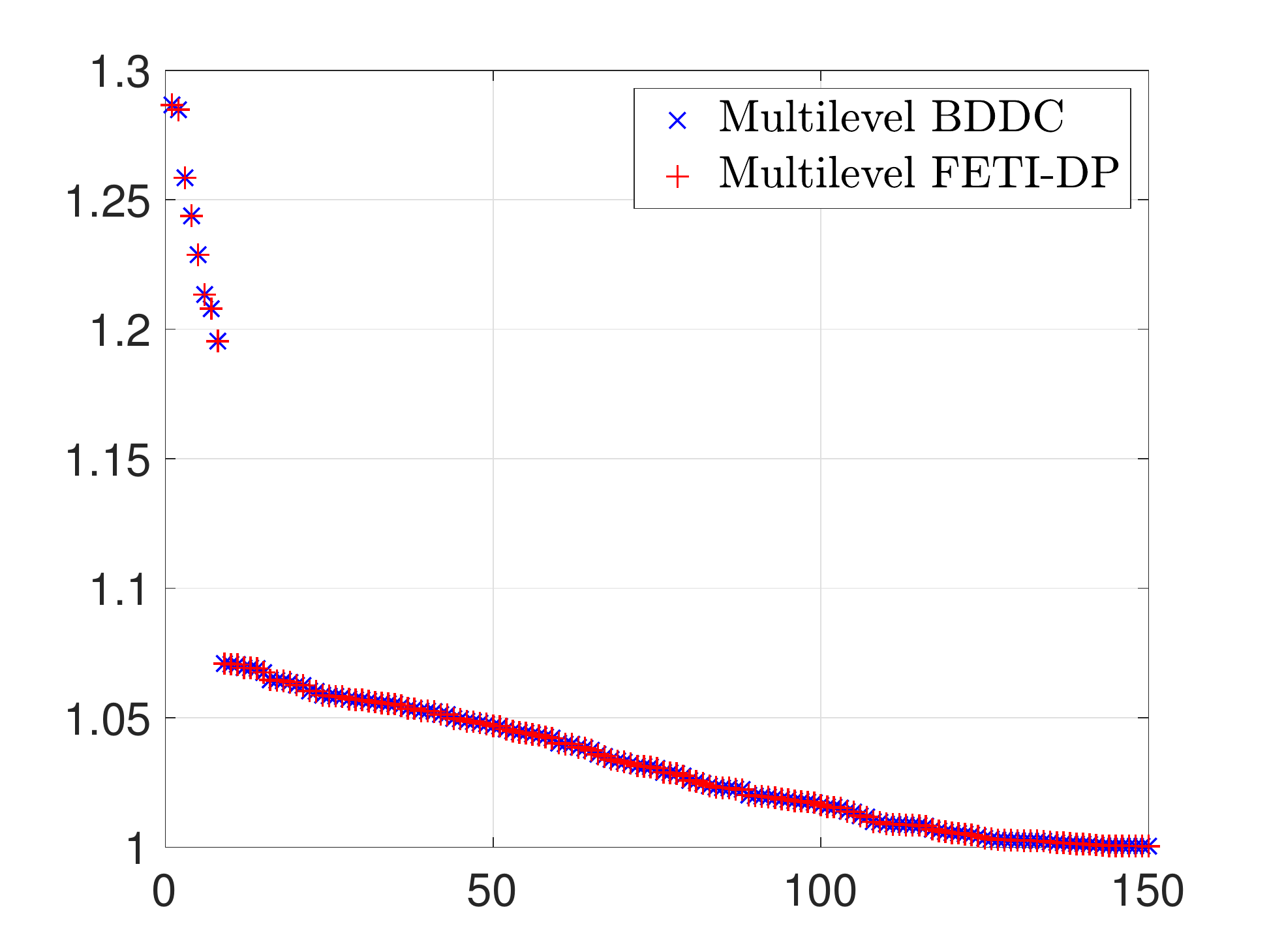}
\end{tabular}
\end{center}
\caption{The largest $150$ eigenvalues of the three-level BDDC and FETI-DP
preconditioned operators with coarsening ratio $H_{\ell}/H_{\ell-1}=3$, where
$\ell=1,2$, with corner constraints (left), and corners and arithmetic
averages over \textcolor{black}{edges} (right).}%
\label{fig}%
\end{figure}

Figure~\ref{fig} displays \textcolor{black}{estimates of} the largest $150$ eigenvalues of the multilevel
FETI-DP\ and BDDC preconditioned operators in the setup with three levels and
coarsening ratio $H_{2}/H_{1}=H_{1}/H_{0}=3$. Specifically, the left panel
corresponds to the second row of Table~\ref{tab:c} and the right panel to the
second row of Table~\ref{tab:c+f}. In both panels all eigenvalues match as
predicted by the theory.

\begin{table}[ptbh]
\caption{Numbers of GMRES iterations (and GMRES/PCG for BDDC) for the
multilevel FETI-DP and BDDC methods, and for the FETI-DP with the block
diagonal preconditioner~(BD) from~(\ref{eq:M_BD}), with increasing number of
levels $L$ and coarsening ratios $H_{\ell}/H_{\ell-1}$, $\ell=1,\dots,L-1$
with corner constraints. Here $nsub$ is the number of subdomains on each
level, $ndof$ is the number of degrees of freedom, and $\lambda_{\max}$ is the
largest eigenvalue of the multilevel BDDC and FETI-DP preconditioned
operators.}%
\label{tab:c}
\begin{center}%
\begin{tabular}
[c]{|c|c|c|c|c|c|c|}\hline
$L$ & $nsub$ & $ndof$ & $\lambda_{\max}$ & BD & FETI-DP & BDDC\\\hline
\multicolumn{7}{|c|}{$H_{\ell}/H_{\ell-1}=3$}\\\hline
2 & 9 & 100 & 1.8781 & 14 & 9 & 10/10\\
3 & 81/9 & 784 & 3.2636 & 23 & 14 & 15/15\\
4 & 729/81/9 & 6724 & 5.3709 & 29 & 20 & 21/21\\
5 & 6561/729/81/9 & 59,536 & 8.8857 & 35 & 27 & 28/29\\\hline
\multicolumn{7}{|c|}{$H_{\ell}/H_{\ell-1}=4$}\\\hline
2 & 16 & 289 & 2.1797 & 10 & 6 & 12/11\\
3 & 256/16 & 4225 & 4.1758 & 29 & 18 & 19/19\\
4 & 4096/256/16 & 66,049 & 7.8472 & 37 & 27 & 28/29\\\hline
\multicolumn{7}{|c|}{$H_{\ell}/H_{\ell-1}=6$}\\\hline
2 & 36 & 1369 & 2.7982 & 24 & 14 & 15/15\\
3 & 1296/36 & 47,089 & 6.0284 & 36 & 24 & 25/26\\\hline
\end{tabular}
\end{center}
\end{table}

\begin{table}[ptbh]
\caption{Numbers of GMRES (and GMRES/PCG for BDDC) iterations for the same problems
as in Table~\ref{tab:c}, but here with corner constraints combined with
averages over edges.}
\begin{center}%
\begin{tabular}
[c]{|c|c|c|c|c|c|c|}\hline
$L$ & $nsub$ & $ndof$ & $\lambda_{\max}$ & BD & FETI-DP & BDDC\\\hline
\multicolumn{7}{|c|}{$H_{\ell}/H_{\ell-1}=3$}\\\hline
2 & 9 & 100 & 1.0550 & 8 & 4 & 4/5\\
3 & 81/9 & 784 & 1.2866 & 11 & 8 & 8/8\\
4 & 729/81/9 & 6724 & 1.6980 & 14 & 11 & 11/11\\
5 & 6561/729/81/9 & 59,536 & 2.1212 & 17 & 14 & 14/14\\\hline
\multicolumn{7}{|c|}{$H_{\ell}/H_{\ell-1}=4$}\\\hline
2 & 16 & 289 & 1.1094 & 7 & 4 & 6/6\\
3 & 256/16 & 4225 & 1.4779 & 13 & 9 & 9/9\\
4 & 4096/256/16 & 66,049 & 1.9393 & 17 & 13 & 13/13\\\hline
\multicolumn{7}{|c|}{$H_{\ell}/H_{\ell-1}=6$}\\\hline
2 & 36 & 1369 & 1.2280 & 13 & 7 & 7/7\\
3 & 1296/36 & 47,089 & 1.7775 & 17 & 11 & 11/11\\\hline
\end{tabular}
\end{center}
\label{tab:c+f}%
\end{table}

\begin{figure}[ptbh]
\begin{center}%
\begin{tabular}
[c]{cc}%
\includegraphics[height=5.5cm]{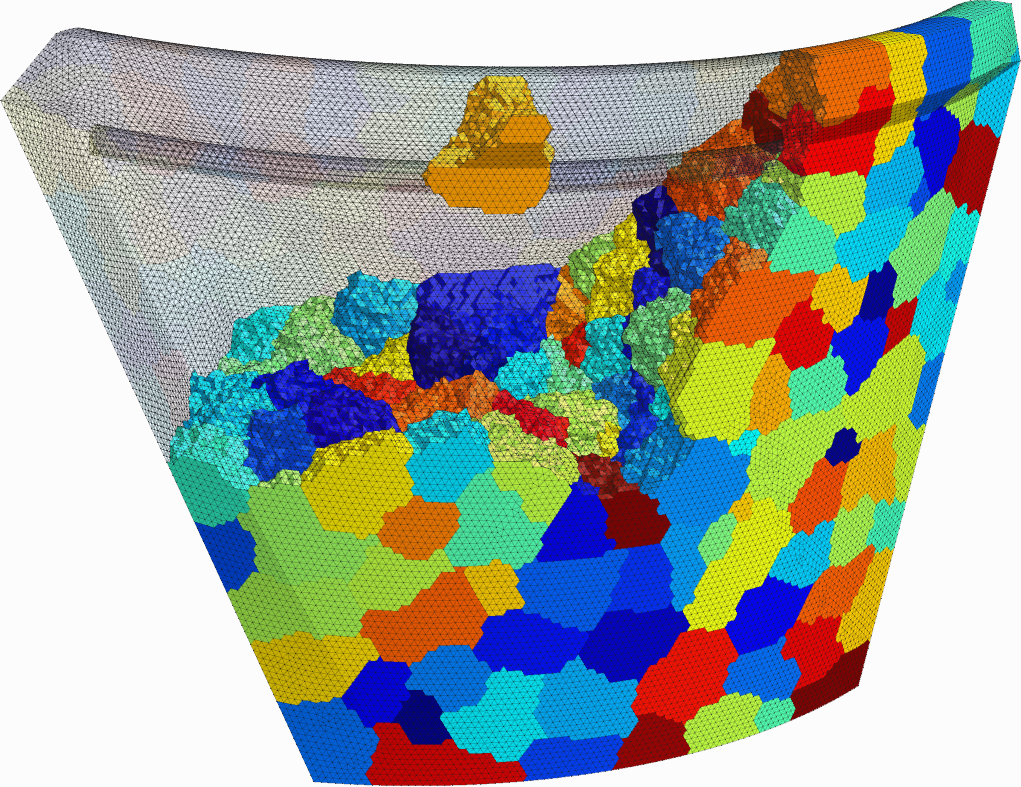} &
\includegraphics[height=5.5cm]{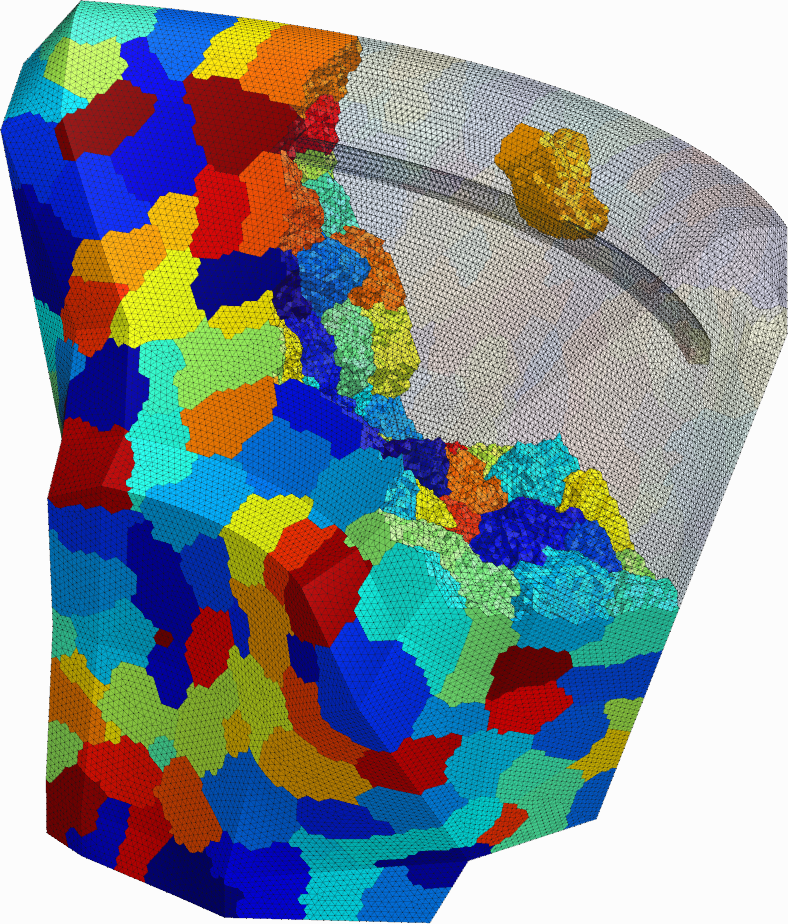}\\
\includegraphics[height=5cm]{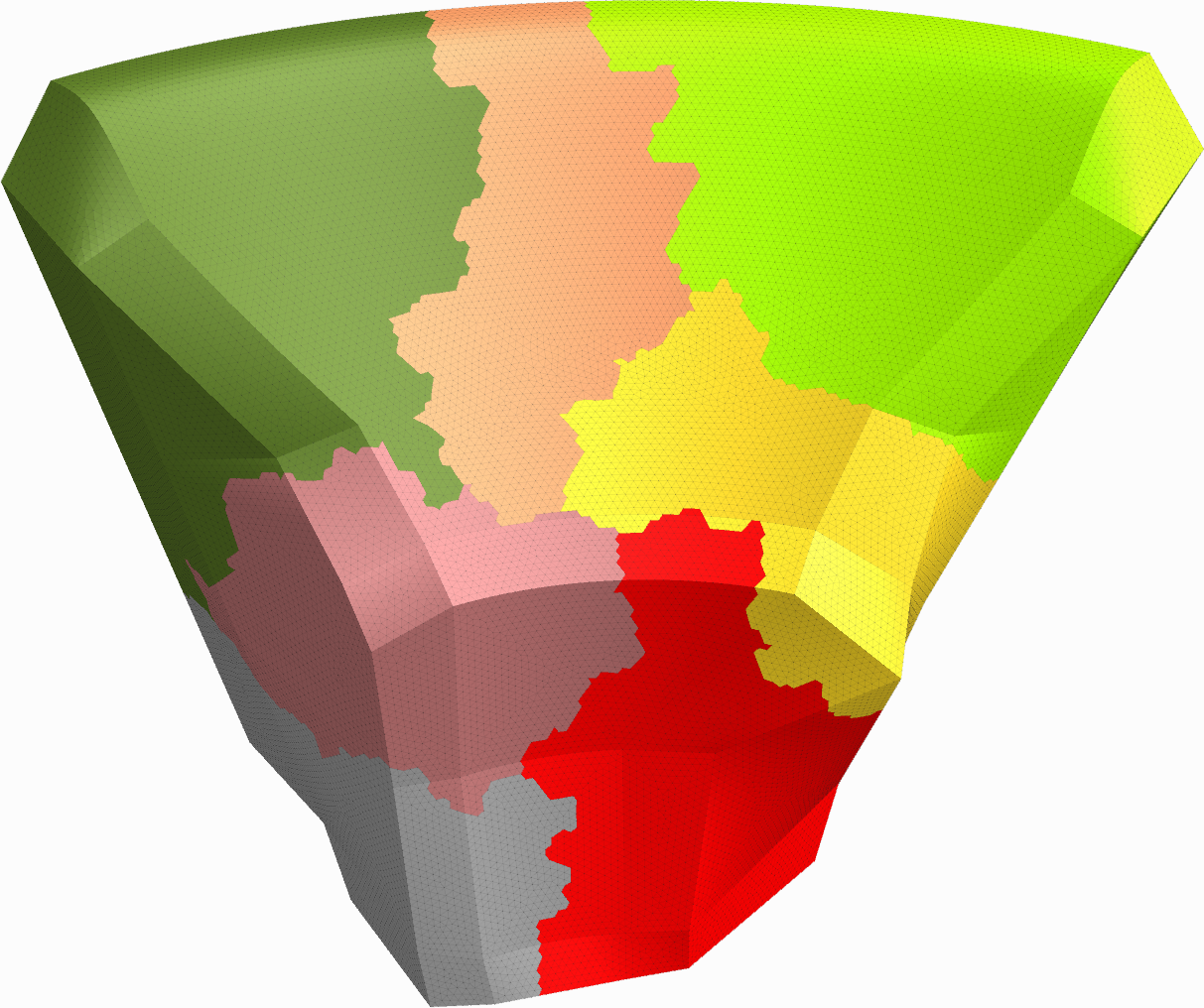} &
\includegraphics[height=5cm]{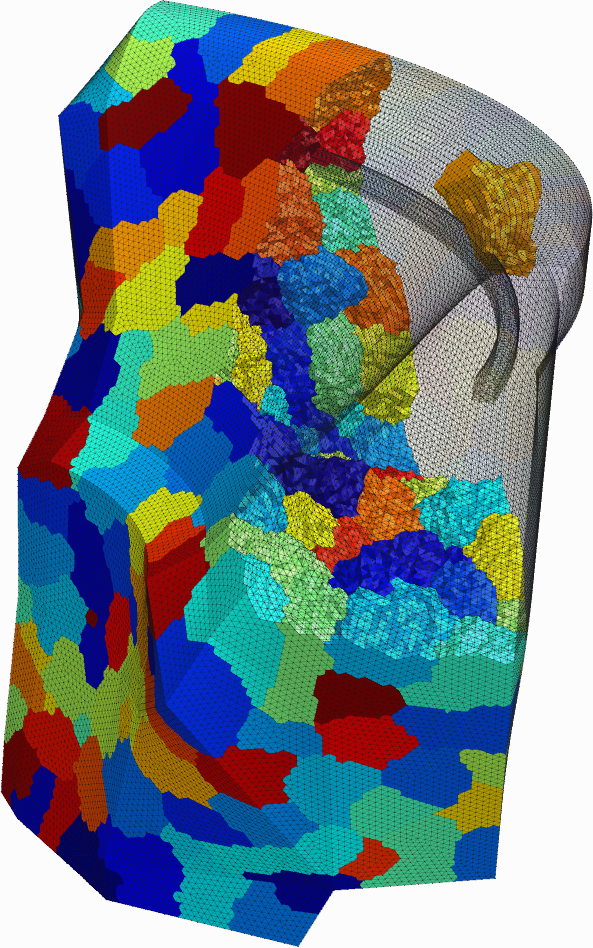}

\end{tabular}
\end{center}
\caption{Finite element discretization and substructuring of the dam
consisting of $3,800,080$ elements distributed into $400$ subdomains in the
first level (top) and $8$ subdomains in the second level (bottom left). Model by
courtesy of Jaroslav Kruis, images reproduced
from~\cite{Sousedik-2010-AMB-thesis}.}%
\label{fig:dam}%
\end{figure}

\begin{figure}[ptbh]
\begin{center}%
\begin{tabular}
[c]{cc}%
\includegraphics[width=5.5cm]{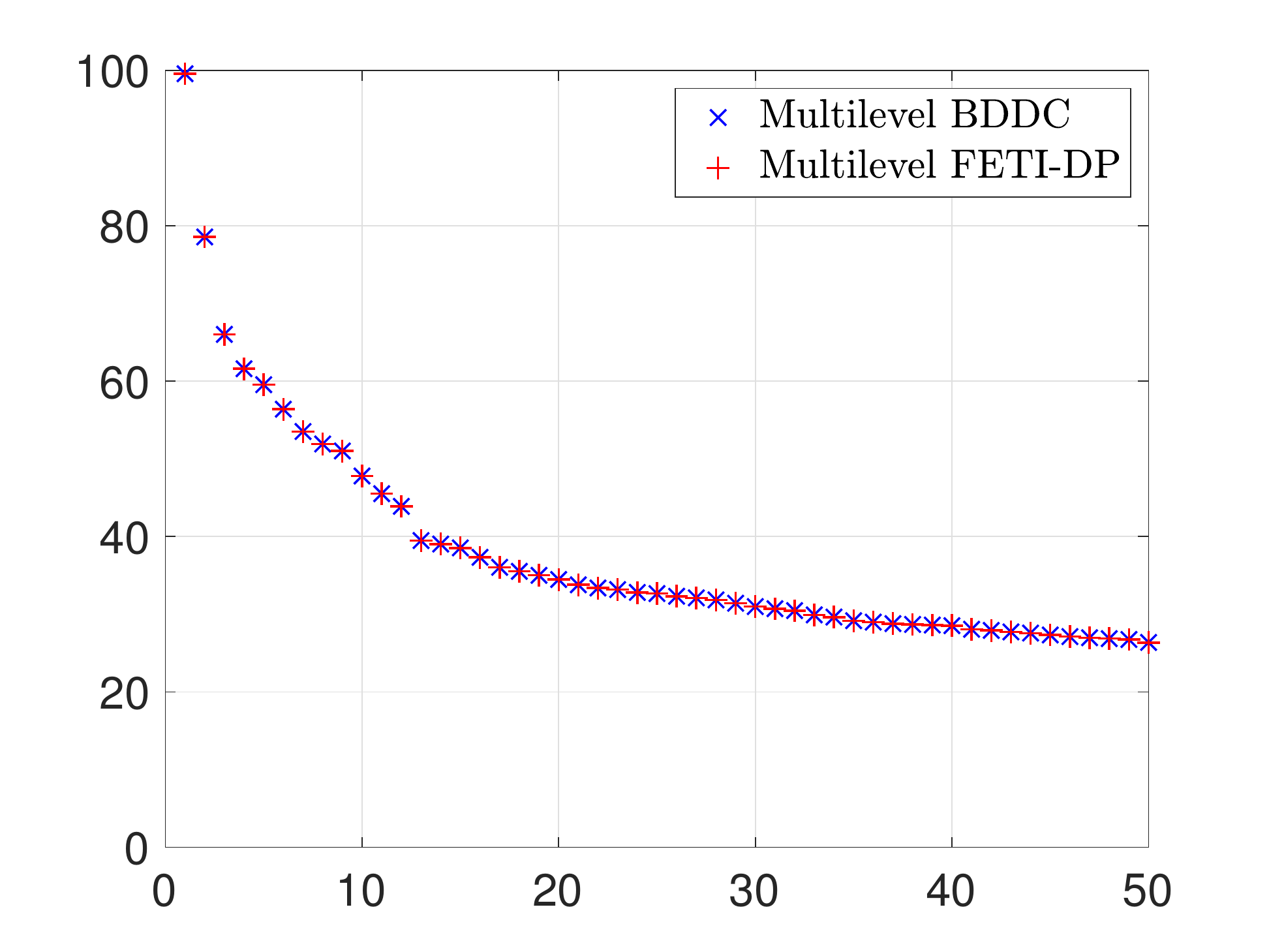} &
\includegraphics[width=5.5cm]{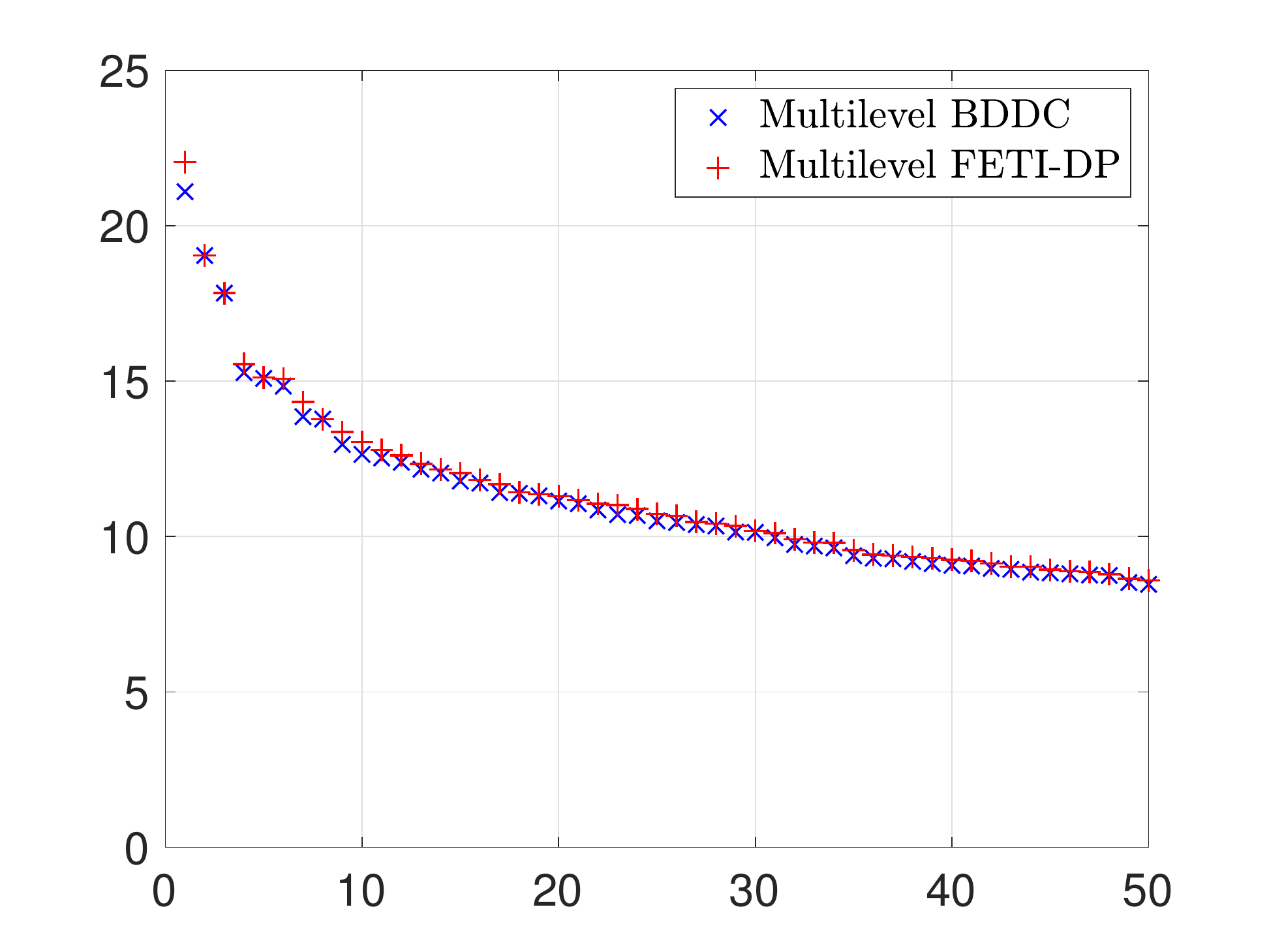}
\end{tabular}
\end{center}
\caption{The largest $50$ eigenvalues of the three-level BDDC and FETI-DP
preconditioned operators for the dam problem using corner constraints (left),
and corners combined with arithmetic averages over edges (right).}
\label{fig:dam-eig}%
\end{figure}

\begin{table}[ptbh]
\caption{Numbers of GMRES iterations (and GMRES/PCG for BDDC) for the dam
problem and the BDDC and FETI-DP methods with varying choices \textcolor{black}{of} constraints:
corners (c), arithmetic averages over edges (e) and faces (f), and their total
numbers on level~1/level~2, and also the largest eigenvalues~$\lambda_{\max}$
of the preconditioned operators. The inexact methods apply algebraic multigrid
(AMG) to the coarse problem. }
\begin{center}%
\begin{tabular}
[c]{|c|c|c|c|c|c|c|c|}\hline
\multicolumn{2}{|c|}{} & \multicolumn{2}{|c|}{inexact methods} &
\multicolumn{4}{|c|}{three-level methods}\\\hline
\multicolumn{2}{|c|}{constraints} & \multicolumn{1}{|c|}{FETI-DP} &
\multicolumn{1}{|c|}{BDDC} & \multicolumn{2}{|c|}{FETI-DP} &
\multicolumn{2}{|c|}{BDDC}\\\hline
type & level~1/level~2 & it & it & $\lambda_{\max}$ & it & $\lambda_{\max}$ &
it\\\hline
c & 11,970/99 & 65 & 71/75 & 99.5 & 66 & 99.5 & 73/87 \\
c+e & 21,180/144 & 105 & 107/125 & 22.0 & 43 & 21.1 & 42/45 \\
c+e+f & 28,002/198 & 134 & 127/140 & 15.1 & 37 & 14.1 & 33/35 \\ \hline
\end{tabular}
\end{center}
\label{tab:dam}%
\end{table}

\textcolor{black}{As the second problem we used a real-world model of linear elasticity in a dam
discretized using $3,800,080$ finite elements with $2,006,748$ degrees of
freedom, distributed into $400$ subdomains with $3990$ corners, $3070$
edges and $2274$ faces. The subdomains in the first level were distributed
into $8$ subdomains for the three-level method.} The distribution of the
elements into subdomains in both substructuring levels was obtained using
\textsc{Metis}~\cite{Karypis-1998-FHQ,Karypis-1998-MSP}, see
Figure~\ref{fig:dam}. Table~\ref{tab:dam} shows the iteration counts of GMRES
(and PCG) for inexact and three-level FETI-DP\ and BDDC\ methods with varying
choices of constraints. The table also shows the largest eigenvalues of the
two preconditioned operators for the three-level methods. For the inexact
methods we applied algebraic multigrid (AMG) as the preconditioner$~M_{c}$
in~(\ref{eq:M_D}) and~(\ref{eq:mBDDC})--(\ref{eq:H-tilde}). Specifically, we
used a \textsc{Matlab} version of the routine \texttt{HSL\_MI20}%
~\cite{Boyle-2010-HSL}, which performs grid coarsening using the classical
algorithm by Ruge and St\"{u}ben~\cite{Ruge-1987-AMG}. We used the
AMG\ routine with its default settings. It can be seen that adding more
constraints \textcolor{black}{as averages over the edges and faces} to the coarse problem (and thus increasing its size) leads to an
increase in the iteration count for the inexact methods. Nevertheless with the
corner constraints, the inexact methods converge in fewer iterations than the
three-level methods. On the other hand, it is expected that adding more constraints would improve convergence for the multilevel methods. 
We see that this \textcolor{black}{happens, in particular after combining the corners with the averages over edges,
and the improvement is less dramatic after adding the averages over faces. 
We also note that there} 
are adaptive techniques for construction of the the coarse spaces,
in\ particular the Adaptive-Multilevel
BDDC~\cite{Mandel-2007-ASF,Sousedik-2013-AMB}, which could be applied to the
multilevel methods presented here, but this is beyond the focus of our study.
Also, a thorough comparison of the proposed methods with AMG\ would be of
independent interest. For example, in a recent study Klawonn et
al.~\cite{Klawonn-2019-PCP} found that, for the BDDC method, inexact
preconditioning\ using AMG\ and the two-level BDDC applied to the coarse
problem show in most situations a very similar behavior. 

\textcolor{black}{Figure~\ref{fig:dam-eig} displays estimates of the $50$ largest eigenvalues for the
three-level FETI-DP\ and BDDC\ preconditioned operators for the setup using
corner constraints (left panel), and using corners combined with arithmetic
averages over edges (right panel), for the dam problem. 
We see that adding averages of the edges reduces the largest eigenvalue approximately four times,
and even though the eigenvalues in the case of the corner constraints decrease at a somewhat faster rate,
their values remain higher than those in the case when both corner and edges are used for constraints. 
In both panels a (close) correspondence of the eigenvalue estimates can be observed for all eigenvalues.}

\begin{figure}[ptbh]
\begin{center}%
\begin{tabular}
[c]{cc}%
\includegraphics[width=5.5cm]{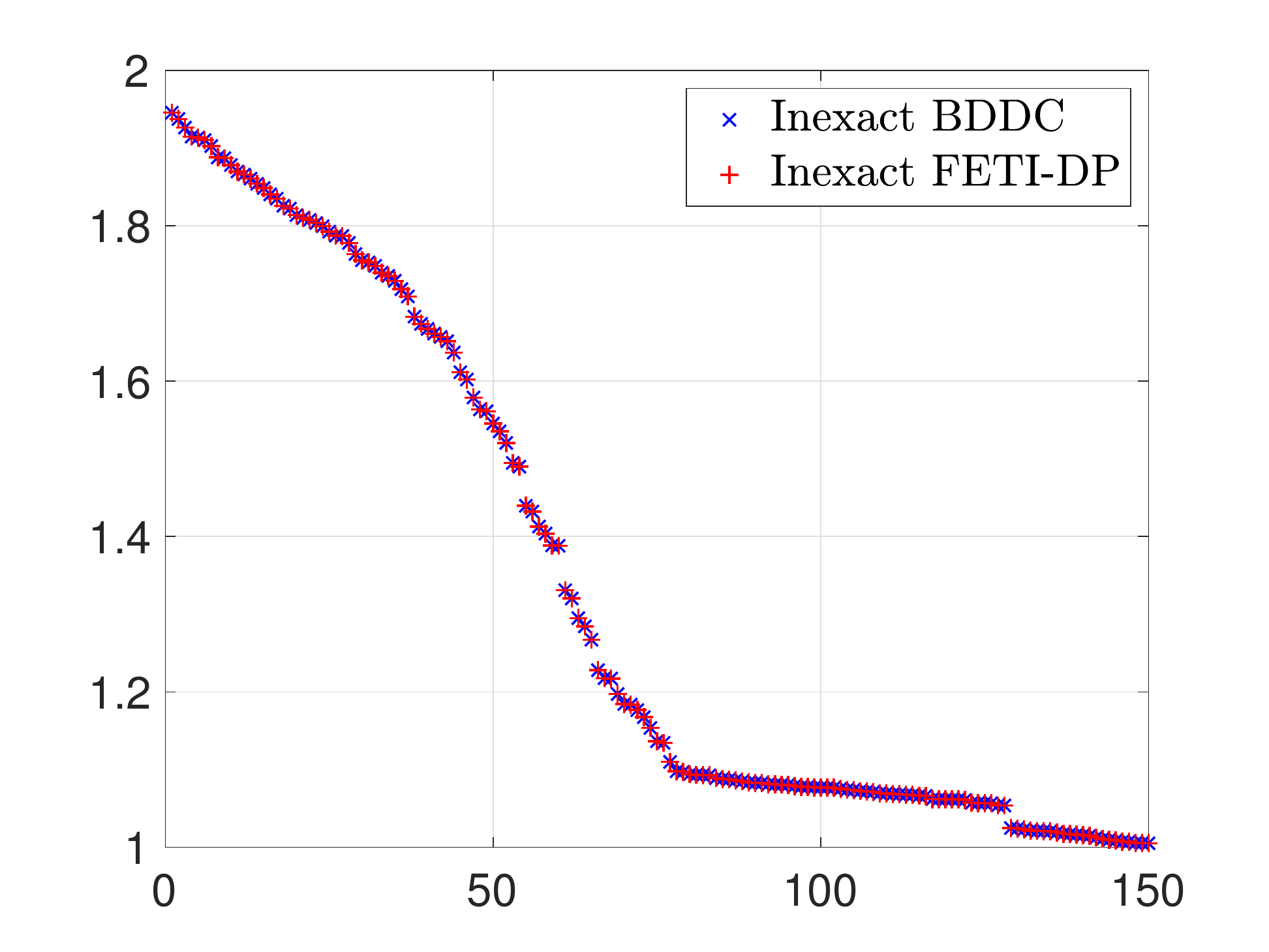} &
\includegraphics[width=5.5cm]{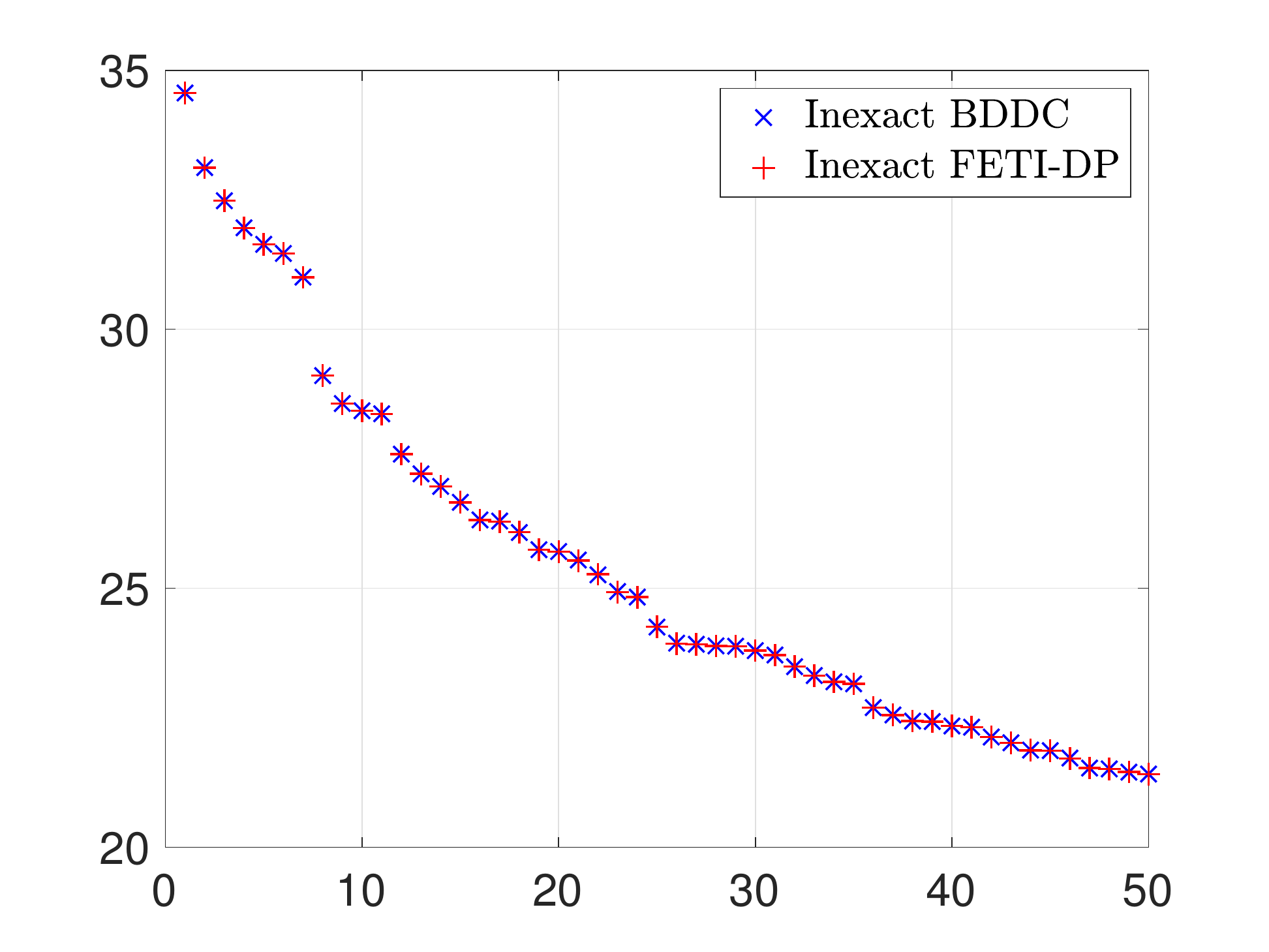}
\end{tabular}
\end{center}
\caption{The largest $150$ eigenvalues of the inexact BDDC and FETI-DP
preconditioned operators with geometric discretization and partitioning of the
domain (left panel), and $50$ eigenvalues of the largest magnitude 
for the dam problem (right).}%
\label{fig:amg}%
\end{figure}

Finally, we note that the iteration counts of GMRES\ and PCG in
Table~\ref{tab:dam}\ are close also for both inexact methods and each set of
constraints. This suggests that the eigenvalues of the inexact methods may be
also similar. A comparison of \textcolor{black}{the eigenvalue estimates} for the inexact methods is provided
by Figure~\ref{fig:amg}. The left panel corresponds to the problem \textcolor{black}{with the same} setup \textcolor{black}{as} in
the second row of Table~\ref{tab:c}, cf. also the left panel in
Figure~\ref{fig}. For this problem, the inexact FETI-DP\ converged in $11$
iterations and the inexact BDDC\ converged in $12$ iterations using both
GMRES\ and PCG. It can be seen that in this case, with geometric
discretization and uniform partitioning of the domain into subdomains, the \textcolor{black}{eigenvalue estimates} are the same (and real). 
\textcolor{black}{An agreement can be observed also in the right panel, which 
displays $50$
eigenvalues of the largest magnitude\ corresponding to the dam problem with
corner constraints. 
Nevertheless, even though all eigenvalues in Figure~\ref{fig:amg} are real,  
we note that the eigenvalues of the inexact FETI-DP\ preconditioned operator are in general complex, 
and a strategy to make them real was provided in~\cite{Klawonn-2007-IFM}.}

\begin{acknowledgements}
The author is grateful to Jaroslav Kruis for providing the discretization of the dam problem, 
to Marian Brezina for visualizations using his own software package~\emph{MyVis},
\textcolor{black}{and to Jan Mandel for continuing support and interest in the domain decomposition methods. 
He would also like to thank the anonymous reviewers for comments and suggestions that helped to improve the quality of the manuscript.}
\end{acknowledgements}

\bibliographystyle{spmpsci}
\bibliography{d-bddc}
%

\appendix

\section{Appendix}

Theorem~\ref{thm:equiv} can be found in~\cite[Theorem 26]{Mandel-2005-ATP},
see also \cite[Theorem 3.9]{Brenner-2007-BFW} and \cite[Theorem 12]%
{Sousedik-2008-CDD}. Since its variant we use here as
Lemma~\ref{lem:equiv-mlevel}, which is virtually the same, we include the
proof for completeness. The proof of Theorem~\ref{thm:equiv} is obtained by
dropping the tilde accent mark. That is, one would write $H$ in place of
$\widetilde{H}$, $F$ in place of $\widetilde{F}$, and $M_{\widehat{S}}$ in
place of $\widetilde{M}_{\widehat{S}}$.

First, let us note a few properties of the operators: both $ER$ and
$B_{D}B^{T}$ are projections,
\begin{align}
BR  &  =0,\quad\text{(from the definition),}\label{eq:BR}\\
REB_{D}^{T}  &  =\left(  I-B_{D}^{T}B\right)  B_{D}^{T}=B_{D}^{T}-B_{D}%
^{T}BB_{D}^{T}=0,\label{eq:REBDt}\\
EB_{D}^{T}B  &  =E\left(  I-RE\right)  =E-ERE=0. \label{eq:EBDtB}%
\end{align}
Next, we show that
\begin{align*}
T_{P}\widetilde{M}_{\widehat{S}}\widehat{S}  &  =M_{D}\widetilde{F}%
T_{P},\qquad T_{P}=B_{D}SR,\\
T_{D}M_{D}\widetilde{F}  &  =\widetilde{M}_{\widehat{S}}\widehat{S}%
T_{D},\qquad T_{D}=E\widetilde{H}B^{T}.
\end{align*}
The first identity is derived as
\begin{align*}
T_{P}\widetilde{M}_{\widehat{S}}\widehat{S}  &  =\left(  B_{D}SR\right)
E\widetilde{H}E^{T}\left(  R^{T}SR\right)  =B_{D}S\left(  I-B_{D}^{T}B\right)
\widetilde{H}\left(  RE\right)  ^{T}SR=\\
&  =B_{D}S\widetilde{H}\left(  RE\right)  ^{T}SR-B_{D}SB_{D}^{T}%
B\widetilde{H}\left(  I-B_{D}^{T}B\right)  ^{T}SR=\\
&  =B_{D}S\widetilde{H}\left(  RE\right)  ^{T}SR-B_{D}SB_{D}^{T}%
B\widetilde{H}SR+B_{D}SB_{D}^{T}B\widetilde{H}B^{T}B_{D}SR=\\
&  =\left(  B_{D}SB_{D}^{T}\right)  B\widetilde{H}B^{T}\left(  B_{D}SR\right)
=M_{D}\widetilde{F}T_{P},
\end{align*}
where some terms vanished, because using~(\ref{eq:HS}), (\ref{eq:REBDt}) and
(\ref{eq:BR})\ we get
\begin{align*}
B_{D}S\widetilde{H}\left(  RE\right)  ^{T}SR  &  =\left(  R^{T}%
SRE\underbrace{\widetilde{H}S}_{=I}B_{D}^{T}\right)  ^{T}=\left(
R^{T}S\underbrace{REB_{D}^{T}}_{=0}\right)  ^{T}=0,\\
B_{D}SB_{D}^{T}B\widetilde{H}SR  &  =B_{D}SB_{D}^{T}B\underbrace{\widetilde{H}%
S}_{=I}R=B_{D}SB_{D}^{T}\underbrace{BR}_{=0}=0.
\end{align*}
The second identity is derived as
\begin{align*}
T_{D}M_{D}\widetilde{F}  &  =\left(  E\widetilde{H}B^{T}\right)  B_{D}%
SB_{D}^{T}\left(  B\widetilde{H}B^{T}\right)  =E\widetilde{H}\left(
I-RE\right)  ^{T}SB_{D}^{T}B\widetilde{H}B^{T}=\\
&  =E\widetilde{H}SB_{D}^{T}B\widetilde{H}B^{T}-E\widetilde{H}\left(
RE\right)  ^{T}S\left(  I-RE\right)  \widetilde{H}B^{T}=\\
&  =E\widetilde{H}SB_{D}^{T}B\widetilde{H}B^{T}-E\widetilde{H}\left(
RE\right)  ^{T}S\widetilde{H}B^{T}+E\widetilde{H}\left(  RE\right)
^{T}SRE\widetilde{H}B^{T}=\\
&  =\left(  E\widetilde{H}E^{T}\right)  R^{T}SR\left(  E\widetilde{H}%
B^{T}\right)  =\widetilde{M}_{\widehat{S}}\widehat{S}T_{D},
\end{align*}
where some terms vanished, because using~(\ref{eq:HS}), (\ref{eq:EBDtB}) and
(\ref{eq:BR}) we get
\begin{align*}
E\widetilde{H}SB_{D}^{T}B\widetilde{H}B^{T}  &  =E\underbrace{\widetilde{H}%
S}_{=I}B_{D}^{T}B\widetilde{H}B^{T}=\underbrace{EB_{D}^{T}B}_{=0}%
\widetilde{H}B^{T}=0,\\
E\widetilde{H}\left(  RE\right)  ^{T}S\widetilde{H}B^{T}  &  =\left(
B\underbrace{\widetilde{H}S}_{=I}RE\widetilde{H}E^{T}\right)  ^{T}=\left(
\underbrace{BR}_{=0}E\widetilde{H}E^{T}\right)  ^{T}=0.
\end{align*}
Finally, let $u_{P}$ be an eigenvector of $\widetilde{M}_{\widehat{S}%
}\widehat{S}$ corresponding to eigenvalue $\lambda_{P}$. Then $T_{P}u_{P}$ is
also an eigenvector of $M_{D}\widetilde{F}$ provided $T_{P}u_{P}\neq0$. So,
let us assume $T_{P}u_{P}=0$. But then,
\begin{align*}
0  &  =T_{D}T_{P}u_{P}=E\widetilde{H}B^{T}B_{D}SRu_{P}=E\widetilde{H}\left(
I-ER\right)  SRu_{P}=\\
&  =E\underbrace{\widetilde{H}S}_{=I}Ru_{P}-E\widetilde{H}ERSRu_{P}%
=ERu_{P}-E\widetilde{H}ERSRu_{P}=ERu_{P}-\widetilde{M}_{\widehat{S}%
}\widehat{S}u_{P},
\end{align*}
but since $ER$ is a projection, $\lambda_{P}$ can be only equal to $0$ or $1$.
Next, let $u_{D}$ be an eigenvector of $M_{D}\widetilde{F}$\ corresponding to
eigenvalue $\lambda_{D}$. Then $T_{D}u_{D}$ is also an eigenvector of
$\widetilde{M}_{\widehat{S}}\widehat{S}$ provided $T_{D}u_{D}\neq0$. So, let
us assume $T_{D}u_{D}=0$. But then,
\begin{align*}
0  &  =T_{P}T_{D}u_{D}=B_{D}SRE\widetilde{H}B^{T}u_{D}=B_{D}S\left(
I-B_{D}^{T}B\right)  \widetilde{H}B^{T}u_{D}=\\
&  =B_{D}\underbrace{S\widetilde{H}}_{=I}B^{T}u_{D}-B_{D}SB_{D}^{T}%
B\widetilde{H}B^{T}u_{D}=B_{D}B^{T}u_{D}-B_{D}SB_{D}^{T}B^{T}\widetilde{H}%
B^{T}u_{D}=\\
&  =B_{D}B^{T}u_{D}-M_{D}\widetilde{F}u_{D},
\end{align*}
but since $B_{D}B^{T}$ is a projection, $\lambda_{D}$ can be only equal to $0$
or $1$.
\end{document}